\documentclass{elsarticle}

\usepackage{amsmath,amsfonts,mathtools}

\usepackage{diagbox}

\usepackage{graphicx}
\usepackage{color}
\usepackage{array}
\usepackage{url}

\usepackage{enumerate}

\usepackage{hyperref}

\newcommand\scv{\mathrm{SCV}}
\newcommand\dt{\mathrm{d}t}

\newcommand\re{\mathrm{Re}}

\graphicspath{{Figures/}}

\newcommand{\Ignore}[1]{}

\usepackage{amsthm}

\newtheorem{theorem}{Theorem}

\usepackage{tikz}
\usepackage{algpseudocode}
\usepackage{algorithm}
\usepackage{pgfplots}

\begin{document}

\sloppy






\title{
Optimized numerical inverse Laplace transformation%
\tnoteref{t1}}

\tnotetext[t1]{This work is partially supported by the OTKA K-138208.}

\author[1]{Ill\'es Horv\'ath}
\ead{horvath.illes.antal@gmail.com}

\author[2]{Andr\'as M\'esz\'aros}
\ead{meszarosa@hit.bme.hu}

\author[1,2]{Mikl\'os Telek}
\ead{telek@hit.bme.hu}

\address[1]{MTA-BME Information Systems Research Group, Budapest, Hungary}
\address[2]{Department of Networked Systems and Services, Technical University of Budapest, Budapest, Hungary}


\begin{abstract}
Among the numerical inverse Laplace transformation (NILT) methods, those that belong to the 
Abate--Whitt framework (AWF) are considered to be the most efficient ones currently. 
It is a characteristic feature of the AWF NILT procedures that they are independent of the transform function and the time point of interest.

In this work we propose an NILT procedure that goes beyond this limitation and optimize the accuracy of the NILT utilizing also the transform function and the time point of interest.

\Ignore{Illes beirna hogy lecsengo fuggvenyekre jo a megkuzelitesunk, Miklos a robosztussagot hangsulyozna}


\end{abstract}

\maketitle

\Ignore{
TODO:\\
-clear up the N vs n notation\\
- add CME-R to the history and the discussion \cite{[Mesz20a]}!!
}

\section{Introduction}\label{s:intro}

Due to the widespread use of Laplace transforms in various scientific fields \cite{LTapplication}, a large number of numerical inverse Laplace transformation (NILT) methods have been developed. Recent surveys are available, e.g., in \cite{Kuhlman2013,WANG201580,iltpaper}.

Among these methods, the most efficient and widely applied ones belong to a subset which is referred to as Abate--Whitt framework (AWF) \cite{AbateUnif}. 
For a given order $N$, each method in the AWF uses a predefined set of $\eta_k$, $\beta_k$  (potentially complex) coefficients independent of the transform function to invert ($h^*(s)$) and the time point of interest ($T$). 
Based on these parameters the AWF NILT procedure is
$$h_N(T)=\sum_{k=0}^{N-1} \frac{\eta_k}{T} ~h^*\left(\frac{\beta_k}{T}\right),$$
where $h^*(s)$ is the Laplace transform function and $h_N(T)$ is the order $N$ approximate of its inverse transform ($h(t)$) at point $T$. Within this framework, various options are available for selecting the $\eta_k$, $\beta_k$ coefficients, in order to obtain efficient NILT methods.
Some methods, e.g., the currently most well known one, the Euler method \cite{Abate2000},
select the $\eta_k$, $\beta_k$ coefficients to closely approximate the 
Bromwich inversion formula \cite[Theorem 24.4]{doetsch}, while some other methods, 
e.g., the recently published CME method \cite{iltpaper}, optimize the \textit{weight function} defined as
$$f_N(t)=\sum_{k=0}^{N-1} \eta_k e^{-\beta_k t},$$
such that it closely approximates the unit impulse function. 

In this paper, we propose a generalization of the AWF such that the NILT method is optimized also for the given transform function to invert ($h^*(s)$) and for the time point of interest ($T$). 
This proposed approach is composed by the following elements: 
\begin{itemize}
\item \textit{a parametric set of AWF methods} \\
The $\eta_k$, $\beta_k$  coefficients depend on a parameter $\theta$, and 
the $\eta_k(\theta)$, $\beta_k(\theta)$ coefficients define an NILT method of the AWF with NILT procedure 
\begin{equation}
h_N(T,\theta)=\sum_{k=0}^{N-1} \frac{\eta_k(\theta)}{T} ~h^*\left(\frac{\beta_k(\theta)}{T}\right).
\label{eq:NILTalpha}
\end{equation}

\item \textit{an error indicator} ($Err(h_N(T,\theta))$) \\
A parameter computed by a numerical procedure that indicates the error of the approximation $h(T)\approx h_N(T,\theta)$ for a given $h^*(s)$ and $T$. 

\item \textit{an optimization method} \\
A method to find the optimal value of the parameter 
\[\hat{\theta}=\arg\min_\theta Err(h_N(T,\theta)).\]
\end{itemize}


The applicability of our proposed NILT approach is limited to the cases when $h(t)$ is real and nonnegative for $t\geq 0$. 
This assumption holds in many practical applications, e.g., when $h(t)$ represents an intrinsically nonnegative physical quantity like a probability or the level of fluid in a container. 
The framework can also be extended to lower bounded $h(t)$ functions with known lower bound 
$m=\min_{\tau\geq 0} h(\tau)$, since in this case $h(t)+m$ with Laplace transform 
$h^*(s)+m/s$ is a nonnegative function.

The rest of the paper is organized as follows. Section \ref{s:ilt} summarizes the basics of NILT with AWF methods. Section \ref{s:gen} presents a parametric set of AWF methods and discusses its behaviour as a function of the parameter. 
Section \ref{s:measures} provides an error indicator of the CME method, while Section \ref{s:opt} proposes a numerical method to optimize the error indicator, and introduces the optimized CME-S method. Section \ref{s:num} analyses the properties of the CME-S method. 
The optimized version of the Euler method, the Euler-S method, is introduced in Section \ref{sec:Euler-S} and analysed in Section \ref{s:num2}, and finally, Section \ref{s:concl} concludes the paper.

\section{Inverse Laplace transformation and the Abate--Whitt framework}
\label{s:ilt}
\subsection{Inverse Laplace transformation}

The Laplace transform of function $h(t)$ is defined as 
\begin{equation}
h^*(s)=\int_{t=0}^{\infty} e^{-st} h(t) \dt,  
\label{eq:LT}
\end{equation}
where $s$ and consequently $h^*(s)$ are potentially complex valued.
The region of convergence for the integral in \eqref{eq:LT} is always of the form $\{s:\re(s)>a\}$ (possibly including some points of the boundary line $\{s:\re(s) = a\}$), or empty ($a=\infty$), or the entire complex plane ($a=-\infty$). The real constant $a$ is referred to as the abscissa of absolute convergence. 

Based on these properties we can summarize the assumptions applied in this paper:
\begin{itemize}
\item[A1)] $h(t)$ is not known, but it is known to be real and nonnegative for any $t\geq 0$.
\item[A2)] The abscissa of absolute convergence of $h^*(s)$, denoted by $a$, is known. 
\item[A3)] The value of $h^*(s)$ is available for all $s$ such that $Re(s)>a$ and we avoid evaluating $h^*(s)$ for $Re(s)\leq a$. 
\end{itemize}

In a wide range of practically important cases, the symbolic inverse Laplace transform of $h^*(s)$ is not available. In these cases, NILT can be applied to find an approximate value of $h$ at point $T$ (i.e., $h(T)$) based on $h^*(s)$. Currently the most efficient NILT methods belong to the AWF.


\subsection{The Abate--Whitt framework}

The NILT methods of the AWF \cite{AbateUnif} 
approximate the $h(t)$ function in point $T$ as 
\begin{align}\label{eq:approx_formula}
h(T)\approx h_N(T) = 
\sum_{k=0}^{N-1} \frac{\eta_k}{T} h^*\left(\frac{\beta_k}{T}\right) \textrm{, }\quad T>0,
\end{align}
where the \textit{nodes} $\beta_k\, (0\leq k\leq N-1)$ and \textit{weights} $\eta_k\, (0\leq k\leq N-1)$ are real or complex numbers that depend on $N$, but not on the transform function $h^*(s)$ or the time point $T$.  Different nodes and weights 
define different NILT methods of the AWF.


We build on the following \textit{integral interpretation} \cite{iltpaper} of the AWF methods which is obtained from \eqref{eq:approx_formula} by substituting \eqref{eq:LT}:
\begin{align}
\nonumber
h_N(T) &
= \sum_{k=0}^{N-1} \frac{\eta_k}{T} h^*\left(\frac{\beta_k}{T}\right)
= \sum_{k=0}^{N-1} \frac{\eta_k}{T} \int_{0}^{\infty} h(t) \cdot  e^{-\beta_k t/T} \dt
\\&= \int_{0}^{\infty} h(t) \cdot \frac{1}{T} f_N(t/T) \dt 
= \int_{0}^{\infty} h(t T) \cdot  f_N(t) \dt,
\label{eq:hnt}
\end{align}
where 
\begin{align}
\label{eq:fnT}
f_N(t) = \sum_{k=0}^{N-1} \eta_k e^{-\beta_k t}.
\end{align}

That is, the result of an AWF NILT procedure according to \eqref{eq:approx_formula}, is equivalent to the final integral in \eqref{eq:hnt}, where $f_N(t)$ is an appropriately selected \textit{weight function}. 
If $f_N(t)$ was the unit impulse function at one (also referred to as Dirac function), then the integral in \eqref{eq:hnt} would result in a perfect Laplace inversion.
The different AWF methods apply different weight functions
as it is exemplified in Figure \ref{fig:main}. 
The $f_N(t)$ functions of the widely applied AWF methods are such that $\int_0^\infty f_N(t)\dt=1$ and $\mathrm{argmax}_t f_N(t)\approx 1$. 
From the numerous AWF methods, we restrict our attention to the most efficient ones, the Euler and the CME methods. 

Contrary to previous works (e.g. \cite{AbateUnif,iltpaper}), here we avoid the simplification due to the complex conjugate $\eta_k$, $\beta_k$ pairs and
we consider both of them in summations like \eqref{eq:approx_formula}
for the ease of notation.
Since both the Euler and the CME methods have one real and $n$ complex conjugate pairs of nodes we have $N=2n+1$, and they represent order $n+1$ NILT methods requiring the evaluation of $h^*(s)$ in $n+1$ points.

\subsubsection*{Euler method (defined only for even $n$)}

The Euler method is an implementation of the Fourier-series method, using Euler summation to accelerate convergence \cite{Abate2000}. 
We define $N=2n+1$ nodes, such that $\beta_ 0$ is real, 
$\beta_k$ has positive imaginary part for $1\leq k\leq n$ 
and negative imaginary part for $n+1\leq k\leq 2n$. 
The $\eta_k$ weights and $\beta_k$ nodes, for $k = 0,1,\dots,2n$ are as follows:
\begin{align}
\label{eq:eulerbeta}
\beta_0&=\alpha, \quad \beta_k = \alpha + \pi i k, \quad \beta_{2n-k+1} = \alpha - \pi i k, \quad 0 < k\leq n,\\
\eta_0&=e^\alpha, \quad \eta_k = (-1)^k e^{\alpha} \xi_k, \quad \eta_{2n-k+1} = \eta_k, \quad 0 < k\leq n,
\label{eq:eulereta}
\end{align}
where $i=\sqrt{-1}$ is the imaginary unit, $\alpha=\frac{n\ln(10)}{6}$ and
\begin{align*}
\xi_k &= 1,\quad 1\leq k\leq n/2,\\
\xi_{n} &= \frac{1}{2^{n/2}},\\
\xi_{n-k} &= \xi_{n-k+1} + 2^{n/2}{n/2 \choose k}, ~~~~~\textrm{ for } 0 < k < n/2.
\end{align*}

The main properties of the Euler method  are as follows:
\begin{itemize}
\item $f_N(t)$ alternates between positive and negative peaks (c.f. Figure \ref{fig:main}).
\item The main part of $f_N(t)$ has significant waves next to the main peak at $t=1$. 
\item $f_N(t)$ is a product of an exponential decay $e^{-\alpha t}$ and a periodic function, whose period is $2$ (c.f. Figure \ref{fig:main}a)). 
\item The initial part of $f_N(t)$ is flat and is close to $0$ (c.f. Figure \ref{fig:main}a)). 
\item $\max_{1\leq k\leq n}|\eta_k|$ increases exponentially with $n$. 
\end{itemize}
\subsubsection*{The CME method}
\label{ss:applilt}

The CME method \cite{iltpaper} is based on the trigonometric -- exponential relation 
\begin{align}
f_{N}(t) &= c\, \mathrm{e}^{-\lambda t} \prod_{j=1}^{n}{\cos^2\left(\frac{\omega \lambda t - \phi_j}{2}\right)}=
\sum_{k=0}^{N-1} \eta_k e^{-\beta_k t},
\label{eq:meodd2}
\end{align}
whose details are provided in \cite{cmepaper}.
In \eqref{eq:meodd2}, $N=2n+1$, $\beta_ 0$ is real, 
$\beta_k$ has positive imaginary part for $1\leq k\leq n$ 
and complex conjugate negative imaginary part for $n+1\leq k\leq 2n$. 

In the CME method, the $\omega$, $\phi_j, j=1,\ldots,n$ parameters are numerically optimized to minimize the squared coefficient of variation (SCV) 
\begin{align}
\scv:=\frac{\int_{t=0}^\infty t^2f_N(t)\dt \int_{t=0}^\infty f_N(t)\dt}{\left(\int_{t=0}^\infty tf_N(t)\dt\right)^2}-1
\label{eq:scv}
\end{align}
(which is independent of $\lambda$ and $c$), 
and the scaling and normalizing constants, $\lambda$ and $c$, are set to ensure $\int_t f_{N}(t) dt=\int_t t f_{N}(t) dt=1$.

\begin{figure*}
\begin{minipage}{0.45\textwidth}
\centering	
\includegraphics[width=0.98\columnwidth]{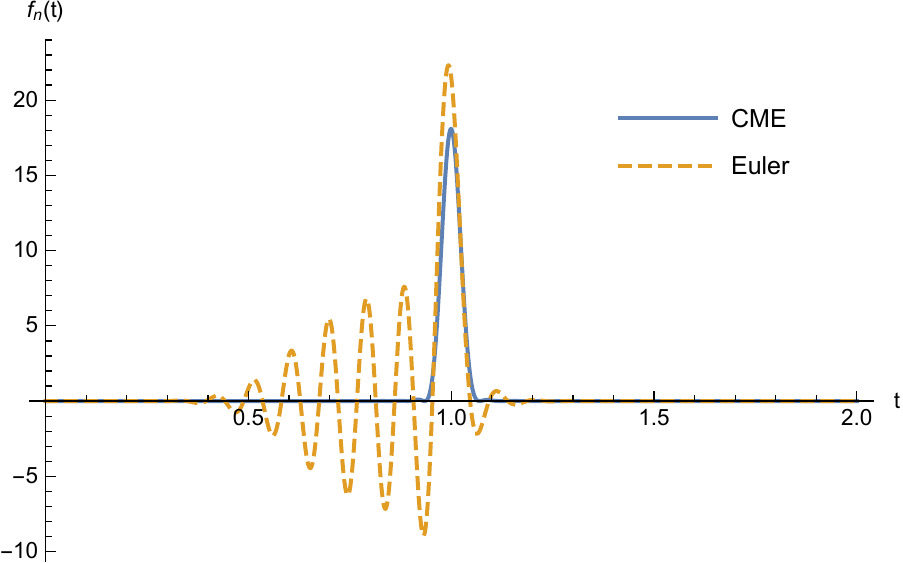}
\\ {\small a) linear y-axis}
\end{minipage}
\hfill
\begin{minipage}{0.45\textwidth}
\centering	
\includegraphics[width=0.98\textwidth]{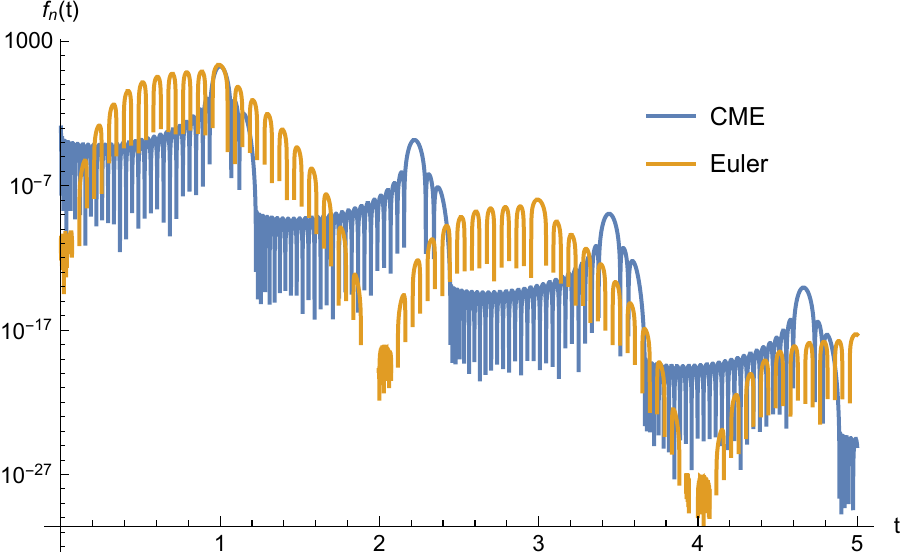}
\\ {\small b) logarithmic y-axis}
\end{minipage}
\caption{The $f_N(t)$ for the Euler and the CME methods for order $30$ with linear  and logarithmic y-axis. The negative parts of the Euler $f_N(t)$ are not visible with logarithmic y-axis.}
\label{fig:main}	
\end{figure*}

The CME method has the following main properties:
\begin{itemize}
\item $f_N(t)$ is nonnegative. 
\item For a given order, the main peak of the CME method is smaller than the one of the Euler method. 
\item The main part of $f_N(t)$ is rather flat apart from the main peak at $t=1$ (c.f. Figure \ref{fig:main}).
\item $f_N(t)$ is a product of an exponential decay and a periodic function (c.f. Figure \ref{fig:main}a)), whose period has no closed form (it is a result of the numerical optimization) and depends on the order. 
\item The initial part of $f_N(t)$ (e.g. between $0$ and $0.1$) has larger peaks than the one of the Euler method (c.f. Figure \ref{fig:main}b)). 
\item $\max_{0\leq k\leq N-1}|\eta_k|$ increases sub-linearly with $N$.
\end{itemize}

\section{A parametric set of Abate--Whitt framework methods}
\label{s:gen}

Let $\eta_k, \beta_k$ be the set of coefficients associated with an AWF method.
For the Euler and the CME methods these coefficients are defined in the previous section 
and for a collection of other AWF methods, they are provided in \cite{AbateUnif}. 
Starting from this set of coefficients, we define
\begin{align}
\eta_k(\theta) =  e^{\theta} \eta_k, ~~~\beta_k(\theta) =   \beta_k + \theta. 
\label{eq:ftheta}
\end{align}
as a function of parameter $\theta$, which we refer to as the shifting parameter. 

To gain an intuitive understanding on the effect of $\theta$ we write the associated weight function as 
\begin{align}
\nonumber
f_{N,\theta}(t)&= \sum_{k=0}^{N-1} \eta_k(\theta) e^{-\beta_k(\theta) t} 
=  \sum_{k=0}^{N-1} (e^{\theta} \eta_k) e^{-(\beta_k+\theta) t} 
\\
\label{eq:fnthetaT}
&=  e^{-\theta (t-1)} \sum_{k=0}^{N-1} \eta_k e^{-\beta_k t} =  e^{-\theta (t-1)} f_N(t). 
\end{align}
Obviously, for $\theta=0$, we obtain the original AWF method with coefficients $\eta_k, \beta_k$. 
If $\theta>0$, then $f_{N,\theta}(t)$ is suppressed for $t>1$ and amplified for $t<1$, compared to $f_{N}(t)$. If $\theta<0$, these relations are reversed.

\begin{figure*}
\begin{minipage}{0.48\textwidth}
\centering	
\includegraphics[width=0.98\textwidth]{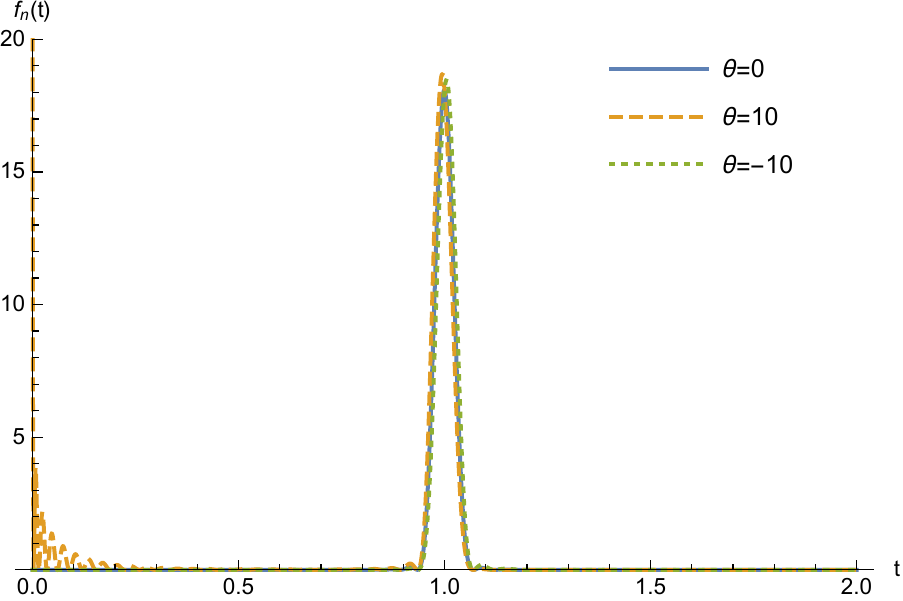}
\\ {\small a) linear y-axis}
\end{minipage}
\hfill
\begin{minipage}{0.48\textwidth}
\centering
\includegraphics[width=0.98\textwidth]{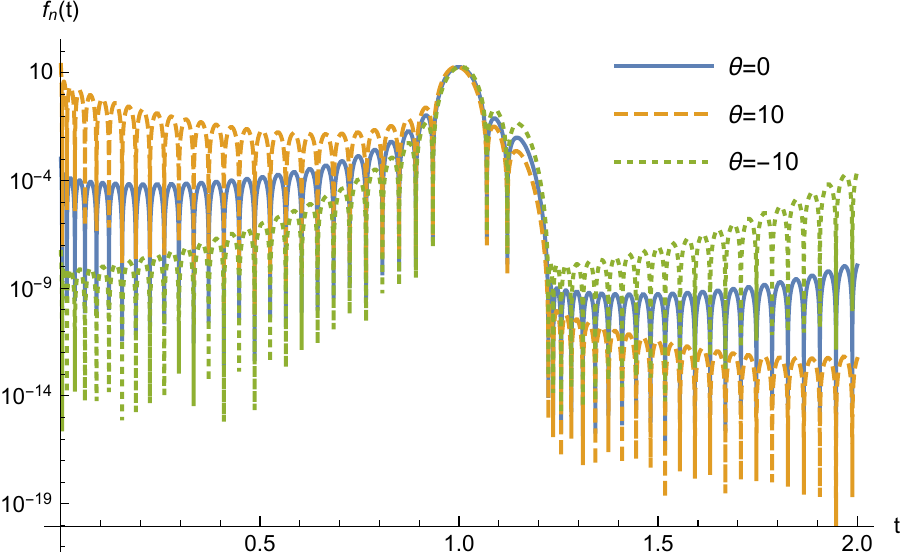}
\\ {\small b) logarithmic y-axis}
\end{minipage}
\caption{The effect of shifting on the weight function of the order $30$ CME method with linear and logarithmic y-axis}
\label{fig:shifting}
\end{figure*}

Figure \ref{fig:shifting} plots the weight functions of the CME method with various shifting parameters with logarithmic and linear y-axis. 
The curves according to logarithmically scaled y-axis in Figure \ref{fig:shifting}b) verifies that a positive shifting parameter amplifies the initial part of the weight function and suppresses its tail, while a negative shifting parameter has the opposite effect. 
Based on Figure \ref{fig:shifting}a), we conclude that the effect of the shifting parameter on the  main peak of the weight function is negligible.  

As a result, the NILT procedure with shifting parameter $\theta$ is as follows
\begin{align}\label{eq:shifting_formula}
h(T)\approx h_N(T,\theta) 
= \sum_{k=0}^{N-1} \frac{\eta_k(\theta)}{T} h^*\left(\frac{\beta_k(\theta)}{T}\right)
=\sum_{k=0}^{N-1} \frac{e^{\theta}\eta_k}{T} h^*\left(\frac{\beta_k+\theta}{T}\right).
\end{align}

The idea of shifting was already introduced in \cite{iltpaper}, where it is recommended that $\theta$ is set equal to the abscissa of convergence $a$, independent of $T$. The approach proposed in this work allows an optimal setting of $\theta$ depending on $h^*(s)$ and $T$.

\section{Error indicator}
\label{s:measures}

In this section, we look for information about the accuracy of an NILT method 
defined by the $\eta_k(\theta), \beta_k(\theta)$ parameters. I.e.,  
the error of the approximation $h(T)\approx h_{N}(T,\theta)$, where $h_{N}(T,\theta)$ is computed according to \eqref{eq:NILTalpha} and $h(T)$ is not known.

\subsection{Properties of the weight functions}

Let $z_1,z_2,\dots$ denote the zeros of $f_{N,\theta}(t)$ for $t > 0$ in increasing order. We set $z_0=0$ (regardless of whether $f_{N,\theta}(0)=0$ or not). According to \eqref{eq:fnthetaT}, the $z_i$ parameters do not depend on $\theta$. 
The index of the largest zero less than one is denoted by $I$, that is, $1\in [z_I,z_{I+1}]$.
We decompose $h_N(T)$ from \eqref{eq:hnt} as
\begin{align}
&h_N(T,\theta)= \int_{0}^{\infty} h(t T) \cdot  f_{N,\theta}(t) \dt 
\label{eq:intparts}= \\
& \underbrace{\int_{0}^{z_I}\!\! h(t T) f_{N,\theta}(t) \dt}_{\varepsilon_\text{left}(\theta)} + \underbrace{\int_{z_I}^{z_{I+1}}\!\! h(t T) f_{N,\theta}(t) \dt}_{h_\text{main}(\theta)} + \underbrace{\int_{z_{I+1}}^{\infty}\!\! h(t T) f_{N,\theta}(t) \dt}_{\varepsilon_\text{right}(\theta)},
\nonumber 
\end{align}
and refer to these terms as the main term, $h_\text{main}(\theta)$, the left error term, $\varepsilon_\text{left}(\theta)$, and the right error term, $\varepsilon_\text{right}(\theta)$.
This naming convention comes from the fact that, if $f_{N,\theta}(t)$ was the unit impulse function at one, then we would have $\varepsilon_\text{left}(\theta)=\varepsilon_\text{right}(\theta)=0$ and $h_\text{main}(\theta)=h_N(T,\theta)=h(T)$.

We can decompose the associated weight functions similarly  
\begin{align}
&\int_{0}^{\infty} f_{N,\theta}(t) \dt =\\
& \underbrace{\int_{0}^{z_I}f_{N,\theta}(t) \dt}_{f_\text{left}} + \underbrace{\int_{z_I}^{z_{I+1}} f_{N,\theta}(t) \dt}_{f_\text{main}} + \underbrace{\int_{z_{I+1}}^{\infty} f_{N,\theta}(t) \dt}_{f_\text{right}} ~~~~(= 1).
\nonumber 
\end{align}

For the Euler weight function $f_\text{main} >> 1$ and $f_\text{left}+f_\text{right} << 0$,
where the $<<$ relation indicates ``significant'' differences. 
In contrast, the CME weight function is nonnegative, consequently, $f_\text{main}$, $f_\text{left}$ and $f_\text{right}$ are all nonnegative, furthermore
$f_\text{main} \approx 1$, 
therefore $1-f_\text{main}<0.01$, as it is demonstrated by Table \ref{tab:prop}. 
For different orders 
$f_\text{main}$, $f_\text{left}$, and $f_\text{right}$ hardly change, while the $(z_I,z_{I+1})$ interval, where the main peak of the weight function is located, decreases significantly with increasing order. 

\begin{table}[h]
\centerline{
\begin{tabular}{|c||c|c||c|c|c|}
\hline
\multicolumn{6}{|c|}{Euler}  \\
\hline
$n$ & $z_I$ & $z_{I+1}$ & $f_\text{left}$& $f_\text{main}$& $f_\text{right}$ \\
\hline
$30$ & 0.9534 & 1.0465 & -0.1492 & 1.1967 & -0.0475 \\
$60$ & 0.9772 & 1.0227 & -0.1528 & 1.2012 & -0.0483 \\
\hline\hline
\multicolumn{6}{|c|}{CME} \\
\hline
$n$ & $z_I$ & $z_{I+1}$ & $f_\text{left}$& $f_\text{main}$& $f_\text{right}$ \\
\hline
$30$ & 0.9344 & 1.0698 & 0.0028 & 0.9950 & 0.0021  \\
$60$ & 0.9689 & 1.0322 & 0.0026 & 0.9949 & 0.0023  \\
\hline
\end{tabular}
}
\caption{Properties of the weight function for the Euler and the CME method
\label{tab:prop}}
\end{table}

According to Assumption A1) $h(t)$ is nonnegative, thus we can interpret the 
$h_\text{main}(\theta)$,  $\varepsilon_\text{left}(\theta)$, $\varepsilon_\text{right}(\theta)$ 
terms depending on the sign of the weight function. 
\begin{itemize}
	\item 
If $f_{N,\theta}(t)$ is nonnegative (like for the CME method and its parametric variants), the terms
$h_\text{main}(\theta)$,  $\varepsilon_\text{left}(\theta)$, and $\varepsilon_\text{right}(\theta)$ are all nonnegative. In this case
$h_\text{main}(\theta)$ approximates $h(T)$, and 
$\varepsilon_\text{left}(\theta)$ and $\varepsilon_\text{right}(\theta)$, represents the error of the  approximation.

\item If  $f_{N,\theta}(t)$ has alternating sign (like in the case of the Euler method and its parametric variants), such clear interpretation of the 
$h_\text{main}(\theta)$,  $\varepsilon_\text{left}(\theta)$, $\varepsilon_\text{right}(\theta)$ terms is not available. 
In this case $h_\text{main}(\theta)>>h(T)$,  $\varepsilon_\text{left}(\theta)<<0$
and $\varepsilon_\text{right}(\theta)<<0$ for ``smooth'' functions (we adopt the intuitively specified concept of smoothness from \cite{Abate2000}).
\end{itemize}

\subsection{Measuring the error by the computed NILT value}
\label{ss:niltvalue}
When both $h(t)$ and $f_{N,\theta}(t)$ are known to be nonnegative, and consequently 
$\varepsilon_\text{left}(\theta)$, $\varepsilon_\text{right}(\theta)$, and $h_\text{main}(\theta)$ are known to be nonnegative, we can approximate the error of the NILT in a computationally efficient way.  

For the parametric Euler and CME methods, 
the main peak of $f_{N,\theta}(t)$ and consequently
$h_\text{main}(\theta)\approx \tilde{h}_\text{main}$ in  \eqref{eq:intparts} is fairly independent of $\theta$, 
as it is exemplified, e.g., in Figure \ref{fig:shifting}.  

For parametric families of AWF methods where the main term in \eqref{eq:intparts} is practically independent of $\theta$ 
\begin{align}
 \min_\theta h_N(T,\theta) &=
\min_\theta 
(\varepsilon_\text{right}(\theta)+ h_\text{main}(\theta)+ \varepsilon_\text{right}(\theta))
\nonumber \\ & \approx
\tilde{h}_\text{main} + \min_\theta (\varepsilon_\text{right}(\theta)+ \varepsilon_\text{right}(\theta)),
\nonumber 
\end{align}
thus minimizing $h_N(T,\theta)$ according to $\theta$ minimizes the error of the NILT as well. 
That is, the NILT value $h_{N,\theta}(T)$ itself can be used to compare the approximation error with different $\theta$ parameters.

\section{Optimization method}
\label{s:opt}

The optimization problem defined in the previous section can be solved with various optimization approaches. To pick a computationally efficient one, we utilize the following 
property of the CME NILT value computed with shifting parameter $\theta$.
\begin{theorem}
\label{th:convex}
If $h(t)$ and $f_{N,\theta}(t)$ are nonnegative functions, then $h_N(T,\theta)$ is a convex function of $\theta$.
\end{theorem}

\begin{proof}
$h_N(T,\theta)$ is convex when $\frac{d^2}{d\theta^2} h_N(T,\theta)\geq 0$.
Substituting the formula for the shifted weight function in \eqref{eq:fnthetaT} into \eqref{eq:hnt} we have 
\begin{align*}
&h_N(T,\theta)= \int_{0}^{\infty} h(t T) \cdot  f_{N,\theta}(t) \dt 
= \int_{0}^{\infty} h(t T) \cdot  e^{-\theta (t-1)} f_{N}(t) \dt
\end{align*}
and
\begin{align*}
& 
\frac{d^2}{d\theta^2} h_N(T,\theta)=  \int_{0}^{\infty} h(t T) \cdot (1-t)^2 e^{\theta (1-t)}  f_{N}(t) \dt \geq 0.
\end{align*}
\end{proof}

That is, to optimize the shifting parameter of the CME based NILT, we have a convex optimization problem to solve.

\subsection{Convex minimization of the computed NILT value}
\label{ss:nilt_impl}

The optimal shifting parameter is obtained as
\begin{align}
\hat{\theta}=\arg\min_\theta h_N(T,\theta),
\label{eq:mintheta}
\end{align}
where $h_N(T,\theta)$ is defined in \eqref{eq:shifting_formula}. 
Consequently, $\hat{\theta}$ is optimized based on $h^*(s)$ and $T$.   
From the solution of \eqref{eq:mintheta}, the proposed NILT approximation is $h_N(T,\hat{\theta})$. 

To find the minimum in \eqref{eq:mintheta}, we make use of the convex behaviour of $h_N(T,\theta)$ in Theorem \ref{th:convex} and apply a simple ternary search method, the 
golden-section search \cite{Kiefer}, where 
the upper and the lower limit of the search method is discussed in the next subsection. 

\subsection{Bounds of the shifting parameter}

In the ternary search optimization method, the initial lower and upper bound for $\theta$ are denoted by $\theta_\ell$ and $\theta_u$, respectively. 
These bounds have to be obtained from $h^*(s)$, $T$ and the original $\eta_k$, $\beta_k$ series (according to \eqref{eq:ftheta}), such that Assumption A3) is met.

If $h^*(s)$ has a finite abscissa of convergence, for the lower bound, we use 
$\theta_\ell=a T-\mu$, where  $a$ is the abscissa of convergence of $h^*(s)$ and $\mu=\max_k Re(\beta_k)$ is the real part of the dominant node of the AWF method.  
This lower bound ensures that all shifted nodes ($\beta_k+\theta$) fall into the convergence region of $h^*(s)$ during the NILT at point $T$, i.e., $\frac{\beta_k+\theta}{T} > a$ for $\forall k$, $\forall \theta>\theta_\ell$.

In the particular case where $a=-\infty$, $\theta_\ell$ is picked arbitrarily, say $\theta_\ell=-1000$; the exact choice is not particularly relevant as long as it is smaller than the optimal $\hat{\theta}$. 
If the ternary search,  starting from this arbitrarily set lower bound, finds the optimal $\hat{\theta}$ to be identical with the lower bound, then
the real optimum  might be lower than the arbitrarily picked lower bound and the ternary search must be restarted from a smaller $\theta_\ell$. 

For the upper bound, we obviously have $\theta_u>aT-\mu$, 
but apart from that it is harder to set. 
If $h(t)$ is known to be bounded, which is the case in many practical applications, $\theta_u=10$ can be used. 
If nothing is known about $h(t)$, then we set $\theta_u$ arbitrarily, say $\theta_u=\max(\theta_\ell+1000,0)$ and if the optimal $\hat{\theta}$ is found to be identical with the arbitrarily set upper bound apply a similar boundary adjustment approach as for the lower bound in case of $a=-\infty$.

\subsection{The proposed NILT procedure}

Putting together the elements from the previous sections, we propose 
Algorithm \ref{alg:cme-s} to enhance of the CME method with shifting, referred to as CME-S, where the optimal shifting parameter is obtained by Algorithm \ref{alg:gss}.

\begin{algorithm}[ht]
\begin{algorithmic}
\Procedure{\textbf{CME-S}}{$h^*(s)$, $T$, $a$, $n$}
\State $\{\boldsymbol{\eta},\boldsymbol{\beta}\}=CMEparams(n)$,
\State $\mu=\mathrm{Re}(\beta_1)$, $\theta_\ell=a T - \mu$, $\theta_h=\max(\theta_\ell+10,10)$, 
\State $\hat{\theta}=\mathbf{GoldenSectionSearch}(h^*(s), T, n, \theta_\ell,\theta_h)$
\State \Return $\mathbf{CME}(h^*(s), T, n, \hat{\theta})=\sum_{k=1}^{N} \frac{e^{\hat{\theta}}\eta_k}{T} h^*\left(\frac{\beta_k+\hat{\theta}}{T}\right)$
\EndProcedure
\end{algorithmic}
\caption{CME method with optimal shifting}
\label{alg:cme-s}
\end{algorithm}

\begin{algorithm}[ht]
\begin{algorithmic}
\Procedure{\textbf{GoldenSectionSearch}}{$h^*(s), T, n, \theta_\ell,\theta_h$}
\State $G=\frac{\sqrt{5}-1}{2}$, $\theta_0=\theta_\ell$, 
 $\theta_1=G \theta_\ell + (1-G) \theta_h$, 
$\theta_2=(1-G) \theta_\ell + G \theta_h$,  $\theta_3=\theta_h$, 
\While {$\theta_3-\theta_0>\varepsilon$}
\If{$\mathbf{CME}(h^*(s), T, n, \theta_1)<\mathbf{CME}(h^*(s), T, n, \theta_2)$}
\State $\theta_3=\theta_2$, $\theta_2=\theta_1$, $\theta_1=G \theta_0 + (1-G) \theta_3$,
\Else
\State $\theta_0=\theta_1$, $\theta_1=\theta_2$, $\theta_2=(1-G) \theta_0 + G \theta_3$,
\EndIf
\EndWhile
\State \Return $\frac{\theta_0+\theta_3}{2}$,
\EndProcedure
\end{algorithmic}
\caption{Golden Section Search method}
\label{alg:gss}
\end{algorithm}

\section{Numerical analysis of the CME-S method}
\label{s:num}

\subsection{Comparing CME and CME-S with regular test functions}

\begin{figure}
\begin{minipage}{0.32\textwidth}
\includegraphics[width=0.98\textwidth]{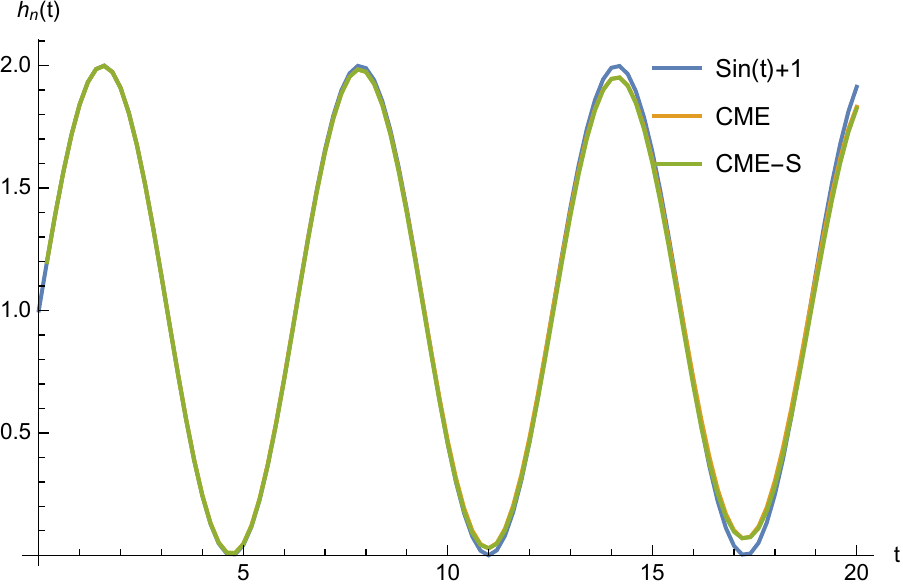}
\end{minipage}
\hfill
\begin{minipage}{0.32\textwidth}
\includegraphics[width=0.98\textwidth]{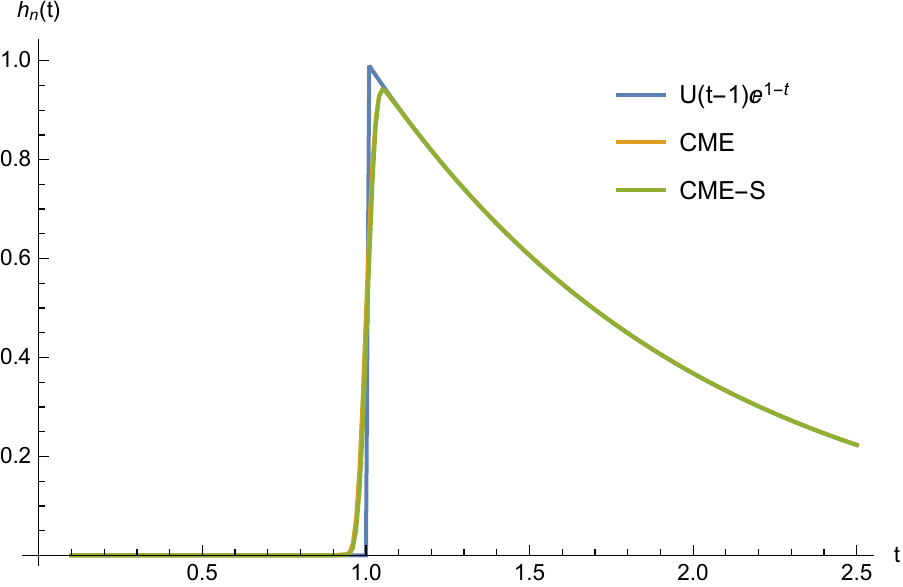}
\end{minipage}
\hfill
\begin{minipage}{0.32\textwidth}
\includegraphics[width=0.98\textwidth]{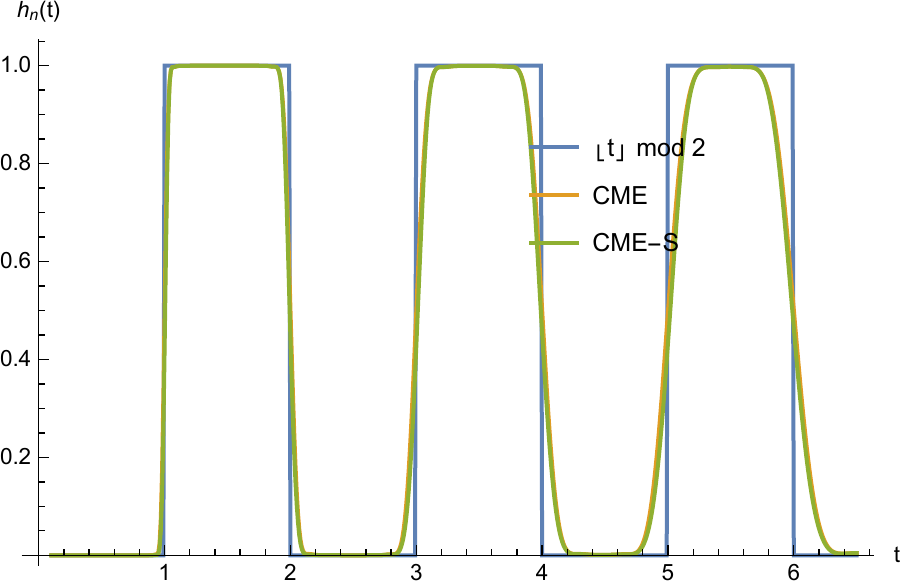}
\end{minipage}
\\
\begin{minipage}{0.32\textwidth}
	\includegraphics[width=0.98\textwidth]{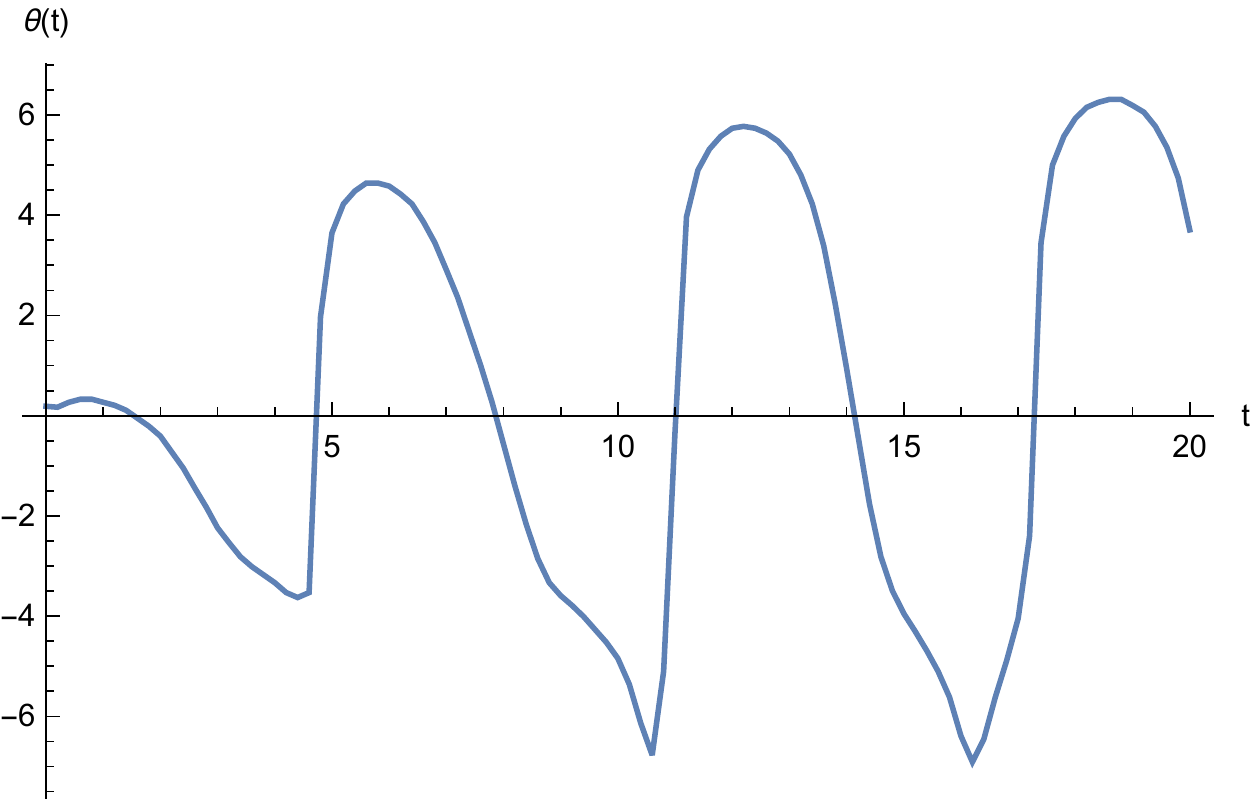}
\end{minipage}
\hfill
\begin{minipage}{0.32\textwidth}
	\includegraphics[width=0.98\textwidth]{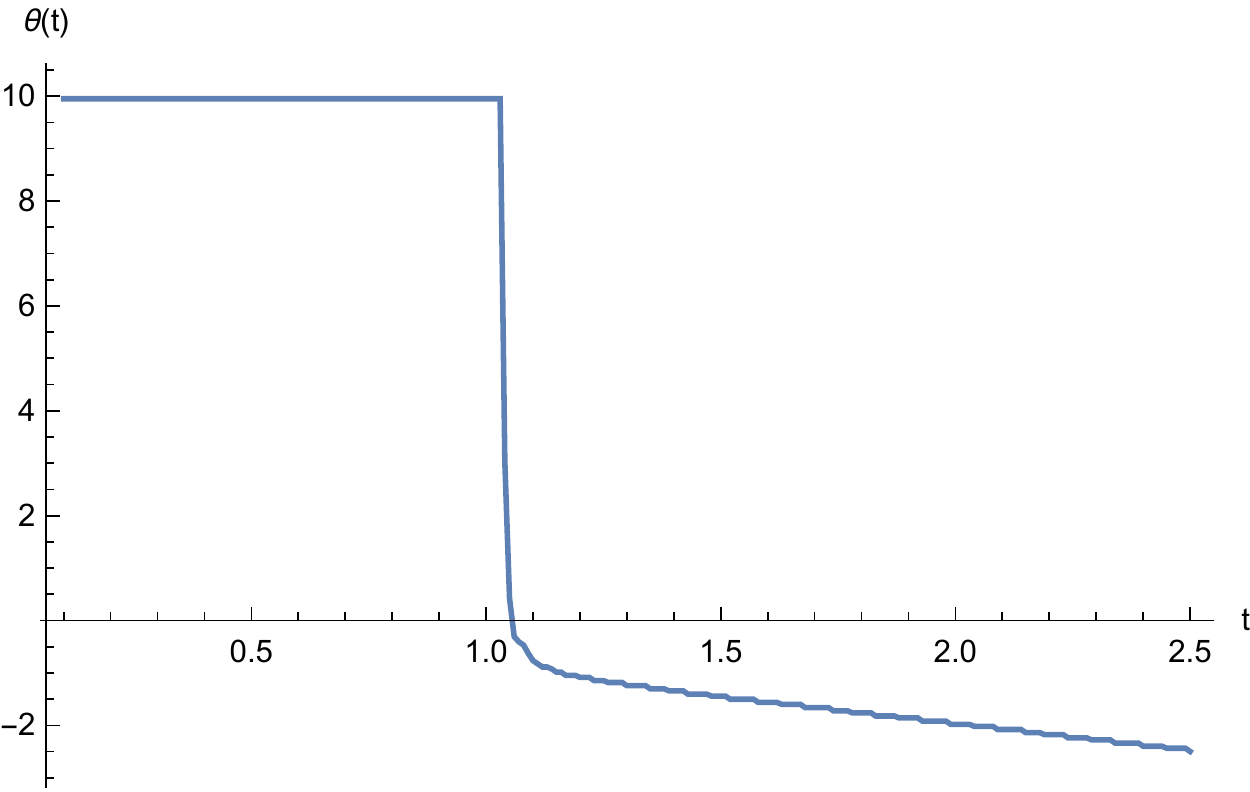}
\end{minipage}
\hfill
\begin{minipage}{0.32\textwidth}
	\includegraphics[width=0.98\textwidth]{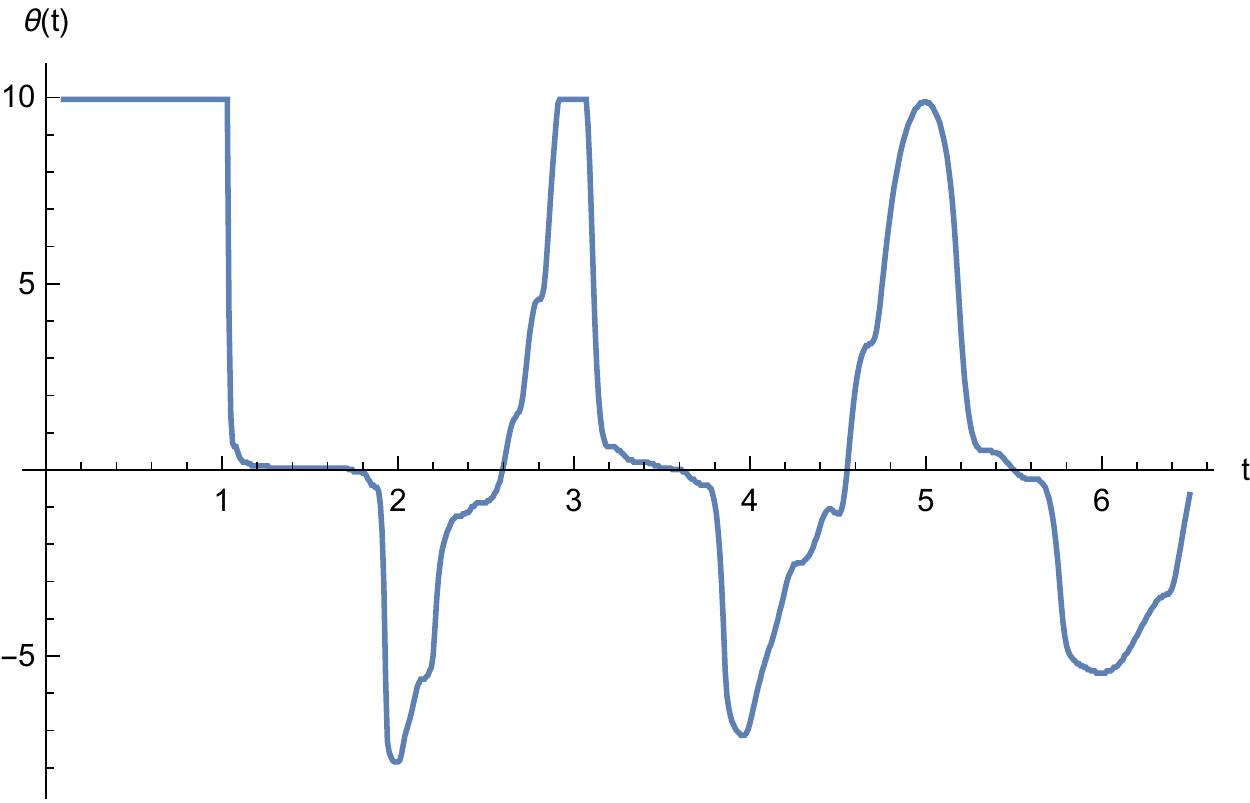}
\end{minipage}
\caption{Behaviour of the CME and the CME-S NILT methods for regular functions $\sin(t)+1$, $U(t-1)e^{1-t}$ and $\lfloor t \rfloor \!\!\mod 2$ with order $30$, where $U(t)$ is the unit step function.}
\label{fig:regfunc}	
\end{figure}

In Figure \ref{fig:regfunc}, we study the behaviour of CME-S for a subset of test functions examined in \cite{iltpaper}: 
$\sin(t)+1 \leftrightarrow \frac{1}{1+s^2}+\frac{1}{s}$
(we use $sin(t)+1$ instead of $sin(t)$ to satisfy the non-negativity of $h(t)$ according to  Assumption A1), 
$U(t-1)e^{1-t}\leftrightarrow \frac{e^{-s}}{1+s}$ and 
$\lfloor t \rfloor \!\!\mod 2 \leftrightarrow \frac{1}{s+s e^s}$. 
Apart from the original function and its CME and CME-S approximation the figure presents the computed $\hat{\theta}$ value of the CME-S procedure. 
In all of these test cases, the abscissa of absolute convergence, $a$, is non-positive and we apply $\theta_h=10$ in the computations.

Based on the results in Figure \ref{fig:regfunc} and several further tests, we conclude that 
\begin{itemize}
\item the CME-S method  does not provide worse results than the CME method, 
\item the $\hat{\theta}$ value computed at point $t$ indicates if $h(t)$ increases or decreases in an environment of point $t$. When $h(t)$ increases around $t$ 
(and consequently the right error is larger then the left error), a positive $\theta$ value helps to decrease the right error, and vice versa. 
\end{itemize}

For a non-negative function which has an initial zero interval (like $h(t)=U(t-1)e^{1-t}$ and $h(t)=\lfloor t \rfloor \!\!\mod 2$ at the $(0,1)$ interval), the optimal shifting would be $\hat{\theta}\to \infty$ in this initial zero interval. 
In this section we have examples where CME-S has similar accuracy as CME, the real benefit of using CME-S is discussed in the next section. 

\subsection{Effect of shifting in tail approximation}

\begin{figure}
\begin{minipage}{0.32\textwidth}
	\includegraphics[width=0.98\textwidth]{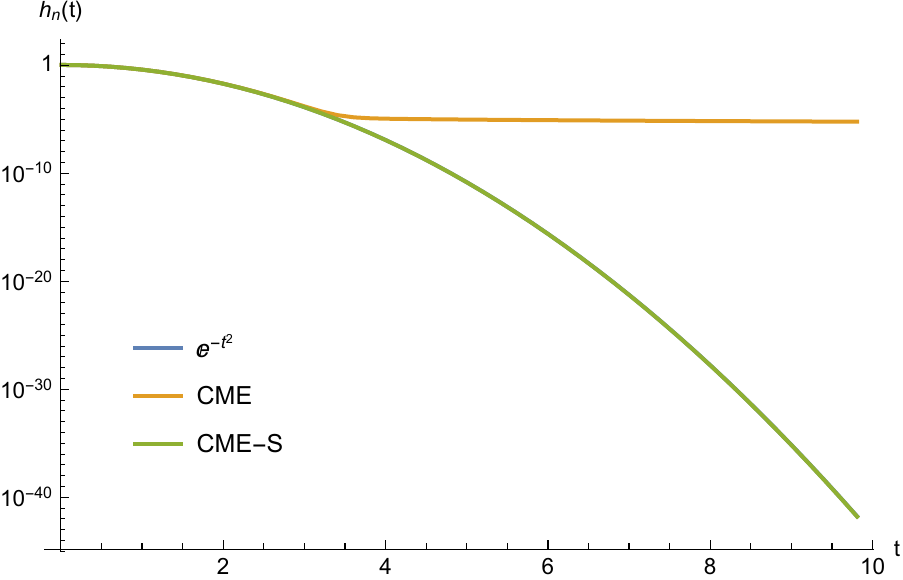}
\end{minipage}
\hfill
\begin{minipage}{0.32\textwidth}
	\includegraphics[width=0.98\textwidth]{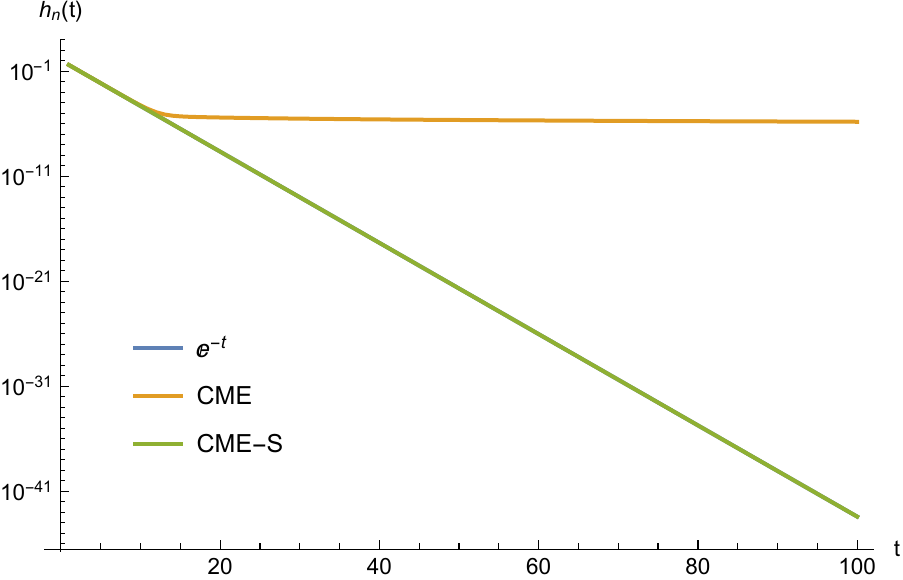}
\end{minipage}
\hfill
\begin{minipage}{0.32\textwidth}
	\includegraphics[width=0.98\textwidth]{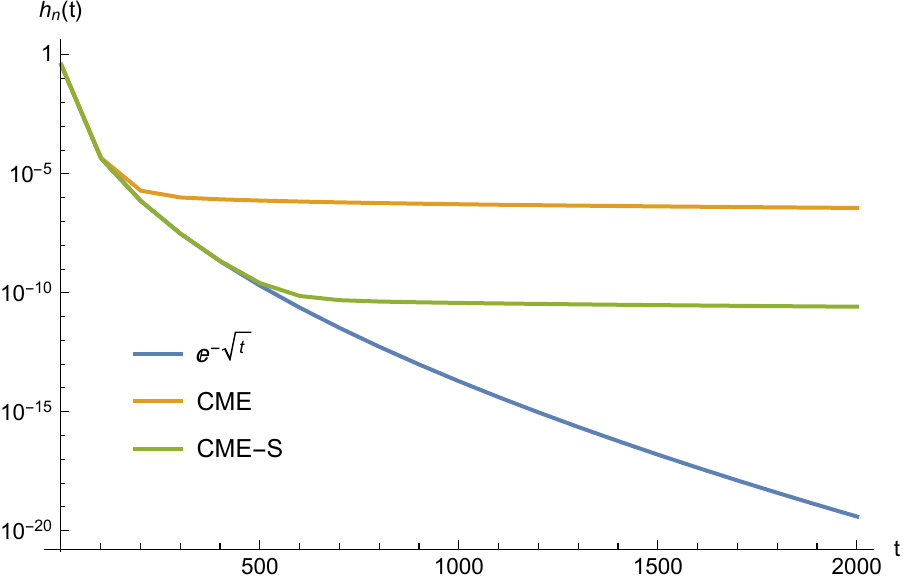}
\end{minipage}
\\
\begin{minipage}{0.32\textwidth}
	\includegraphics[width=0.98\textwidth]{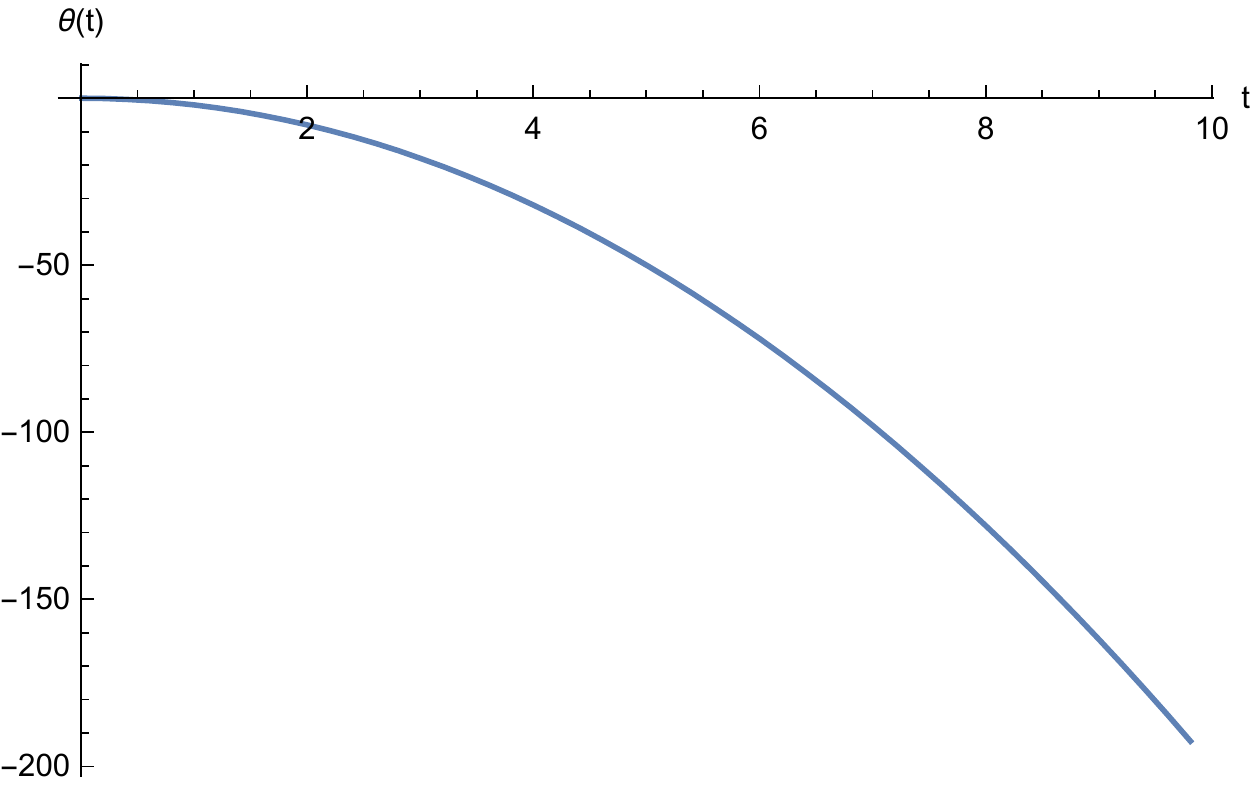}
\end{minipage}
\hfill
\begin{minipage}{0.32\textwidth}
	\includegraphics[width=0.98\textwidth]{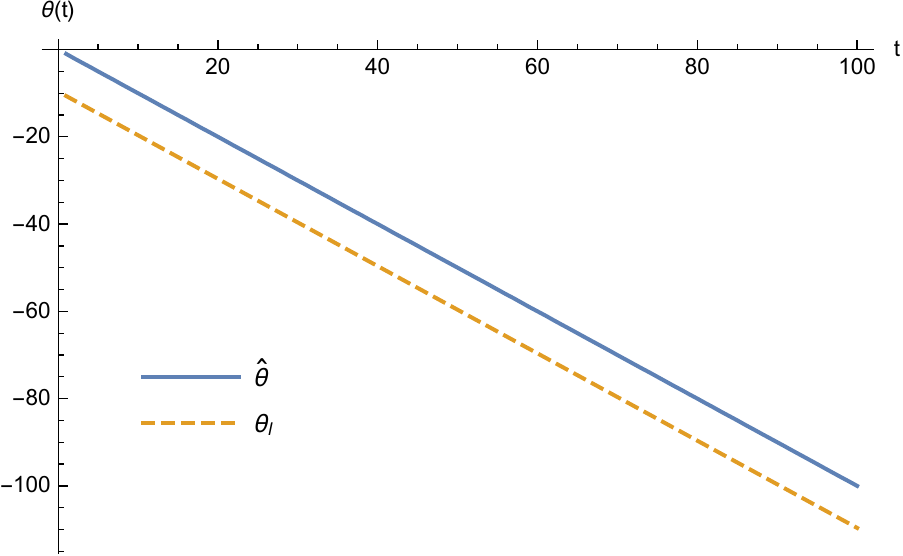}
\end{minipage}
\hfill
\begin{minipage}{0.32\textwidth}
	\includegraphics[width=0.98\textwidth]{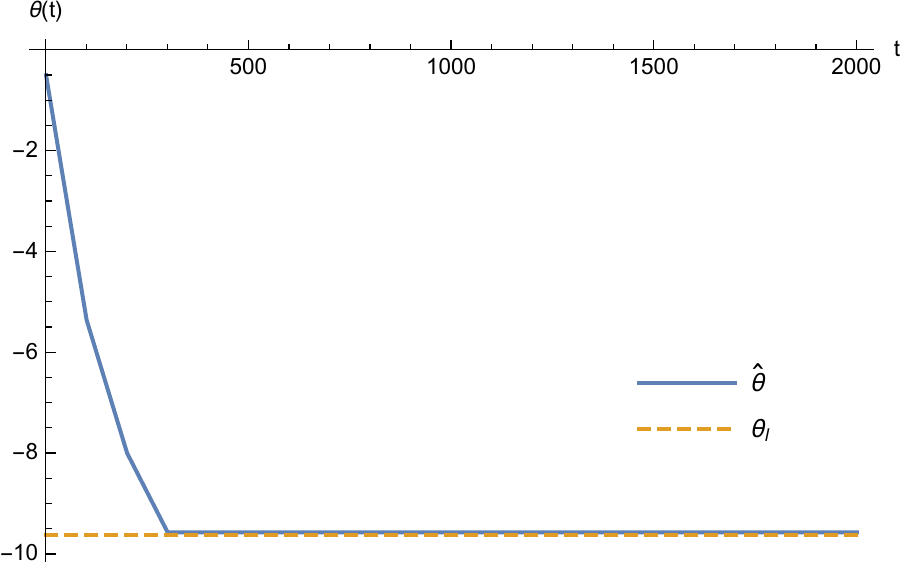}
\end{minipage}
\caption{Behaviour of the CME and the CME-S NILT methods for decaying functions $e^{-t^2}$ ($a=-\infty$), $e^{-t}$ ($a=-1$) and $e^{-\sqrt{t}}$ ($a=0$) with order $30$}
\label{fig:decayfunc}	
\end{figure}

In Figure \ref{fig:decayfunc}, we study the behaviour of CME-S for decaying functions 
in Table \ref{tab:test_functions} for ``large'' $T$ values. In each studied case, optimizing the shifting parameter extends the time interval where the NILT provides correct result. As long as the optimal shifting parameter can follow the decay tendency of $h(t)$ (as it is the case for $e^{-t^2}$ and $e^{-t}$) the CME-S method gives accurate result. When the optimal shifting parameter cannot follow the decay tendency due to the limitation from Assumption A3) (as it is the case for $e^{-\sqrt{t}}$, where the dashed line indicates the $\theta_\ell$ limit) the CME-S method also fails to follow the decay of the original function.

\begin{figure*}
\begin{minipage}{0.49\textwidth}
\centering\includegraphics[width=0.98\textwidth]{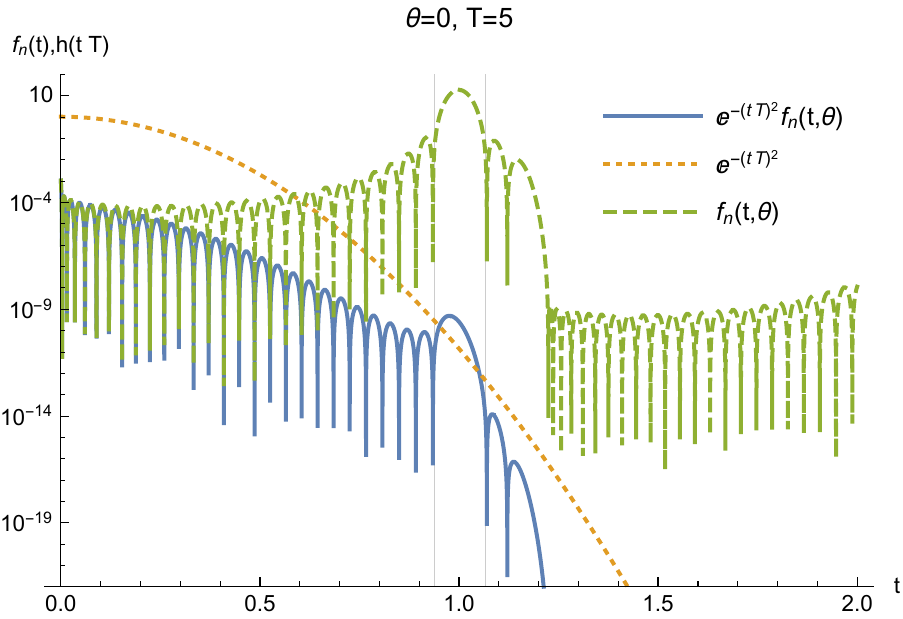}
\\ {\small a) $\theta=0$}
\end{minipage}
\hfill
\begin{minipage}{0.49\textwidth}
\centering\includegraphics[width=0.98\textwidth]{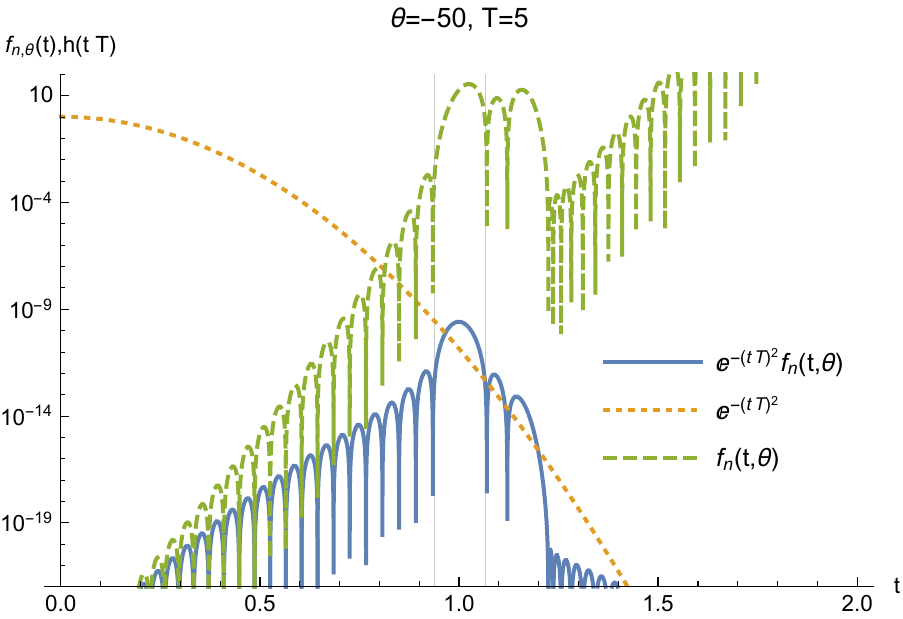}
\\ {\small b) $\theta=-50$}
\end{minipage}
\caption{The elements of the integral interpretation  $h_N(T)=\int_{t=0}^{\infty} h(t T) f_N(t,\theta) dt$ with $h(t)=e^{-t^2}$, $T=5$ and $\theta=0$ and $\theta=-50$ for the CME method with order $30$}
\label{fig:cmeint}
\end{figure*}

The ingredients of the integral interpretation, defined in \eqref{eq:hnt}, 
are depicted in Figure \ref{fig:cmeint} for $h(t)=e^{-t^2}$ and $T=5$.
The figure demonstrates the difficulty of NILT of decaying functions. 
In the plots, the grid lines indicate the integration limits of the left error term, the main term, and the right error term according to \eqref{eq:intparts}. The very sharp decay of 
$h(tT)=e^{-(tT)^2}$ in Figure \ref{fig:cmeint}, 
makes the main term negligibly small compared to the left error term without shifting.
This is why non-optimized NILT methods give many orders of magnitude larger NILT estimates for decaying functions. In Figure \ref{fig:cmeint}b), we shift $f_n(t)$ with $\theta=-50$. In this case $f_{n,\theta}(t)$ is suppressed for small $t$ values and the left error term decreases significantly. The main part hardly changes and the right error increases compared to the non-shifted case. The optimal shifting parameter is the one which makes both, the left and the right errors small compared to the main part.   

Figure \ref{fig:exp2theta}a) plots the computed NILT value as a function of the shifting parameter $\theta$ with the CME method for the same example ($h(t)=e^{-t^2}$ and $T=5$). The figure verifies the convex behaviour, proved in Theorem \ref{th:convex}, for this example.

\section{Euler method with shifting}
\label{sec:Euler-S}

The comparison of the Euler and the CME methods in \cite{iltpaper} indicated that neither of these methods is more accurate than the other in all cases. Based on our qualitative understanding the Euler method is more accurate for ``smooth'' functions, while discontinuities and ``sharp'' changes are better approximated by the CME method for  ``small'' $T$ values, and both of these methods are inaccurate for tail approximation.  

\begin{figure*}
\begin{minipage}{0.49\textwidth}
\centering
\includegraphics[width=0.98\textwidth]{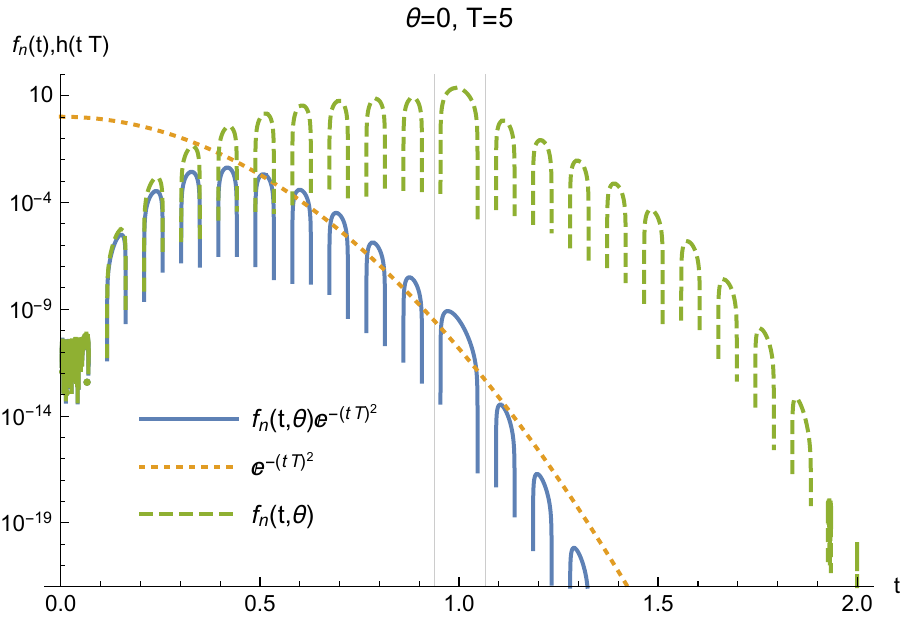}
\\ {\small a) $\theta=0$}
\end{minipage}
\hfill
\begin{minipage}{0.49\textwidth}
\centering
\includegraphics[width=0.98\textwidth]{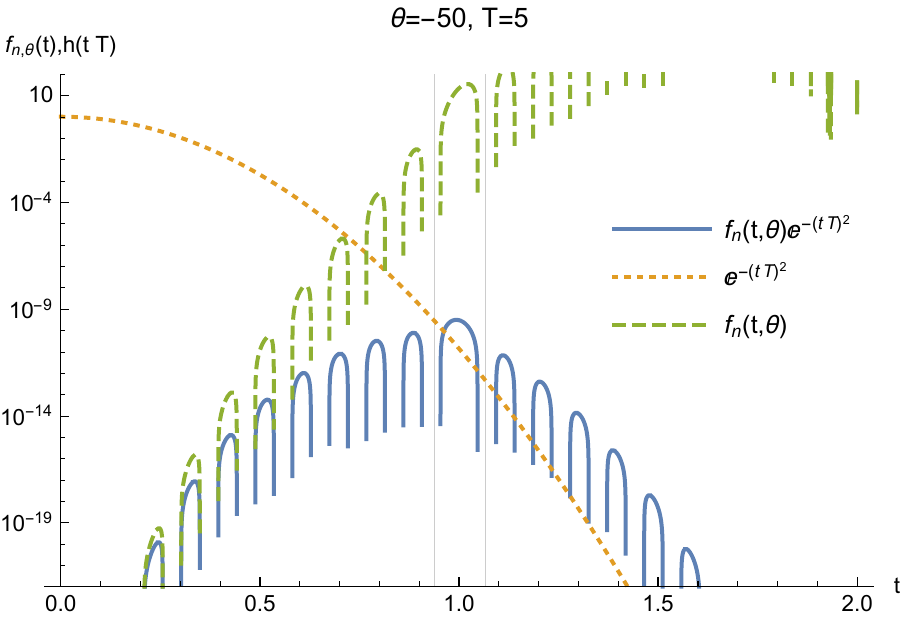}
\\ {\small b) $\theta=-50$}
\end{minipage}
\caption{The elements of integral interpretation  $h_N(T)=\int_{t=0}^{\infty} h(t T) f_N(t,\theta) dt$ with $h(t)=e^{-t^2}$, $T=5$ and $\theta=0$  and $\theta=-50$ for the Euler method with order $30$}
\label{fig:eulerint}
\end{figure*}

Figure \ref{fig:eulerint} demonstrates that 
the Euler method suffers from the same difficulty as the CME method
(c.f. Figure \ref{fig:cmeint}) when it is applied to approximate a decaying function.
The logarithmic scaling of the y-axis makes the negative values invisible in the figure.

\begin{figure*}
\begin{minipage}{0.49\textwidth}
\centering
\includegraphics[width=0.98\textwidth]{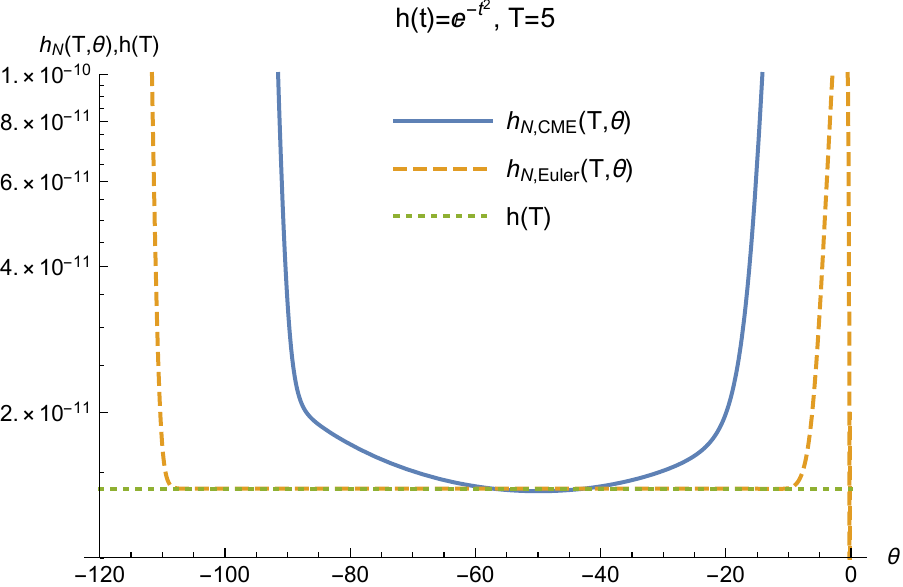}
\\ {\small a) $T=5$}
\end{minipage}
\hfill
\begin{minipage}{0.49\textwidth}
\centering
\includegraphics[width=0.98\textwidth]{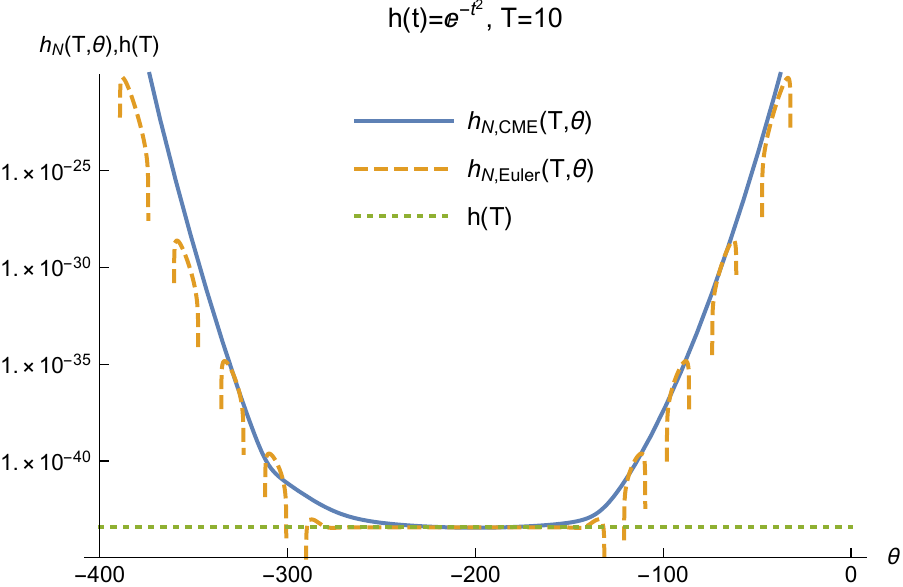}
\\ {\small b) $T=10$}
\end{minipage}
\caption{CME and Euler NILT with various shifting parameters when $h(t)=e^{-t^2}$, $T=5$ and $T=10$ with order $30$}
\label{fig:exp2theta}
\end{figure*}

Figure \ref{fig:exp2theta}a) plots the computed NILT value for 
$h(t)=e^{-t^2}$ and $T=5$ as function of the shifting parameter also with the Euler method. The plot indicates the following properties: 
\begin{itemize}
\item The computed NILT value is a non-convex function of the shifting parameter which might have alternating sign (e.g., it is negative at $\theta=0$ in  Figure \ref{fig:exp2theta}a)). 
\item There is a wider range of $\theta$ values for which the Euler method is reasonably accurate (i.e., a wider range than in case of the CME method). 
\item The range of $\theta$ values where the Euler method is reasonably accurate might be far from zero (the original Euler method is equivalent with $\theta=0$). 
\item The $\theta$ values where the shifted Euler method provide accurate results (i.e., 
$\theta\in(-105,-15)$ in Figure \ref{fig:exp2theta}a) and $\theta\in(-260,-150)$ in Figure \ref{fig:exp2theta}b)) do not have extremal property. 
\item The alternating sign of the computed NILT value makes it hard to find the optimal range of the shifting parameter based on $h_N(t,\theta)$ as a function $\theta$ computed with  Euler method. 
\item The  optimal $\theta$ value of the CME method is (in this example and in Figure \ref{fig:exp2theta}b)) in the $\theta$ range where the Euler method is accurate.   
\end{itemize}

Figure \ref{fig:cmeint} and Figure \ref{fig:eulerint} suggests that 
the optimal $\hat{\theta}$ parameter of the CME method and the $\theta$ range where the Euler method with shifting is accurate coincidence for decaying functions.
As a result, we propose to apply the Euler method with the optimal shifting parameter computed with the CME method. We refer to this extension of the Euler method with shifting as the Euler-S method.  

\begin{figure*}
\hfill	
\begin{minipage}{0.49\textwidth}
	\includegraphics[width=0.98\textwidth]{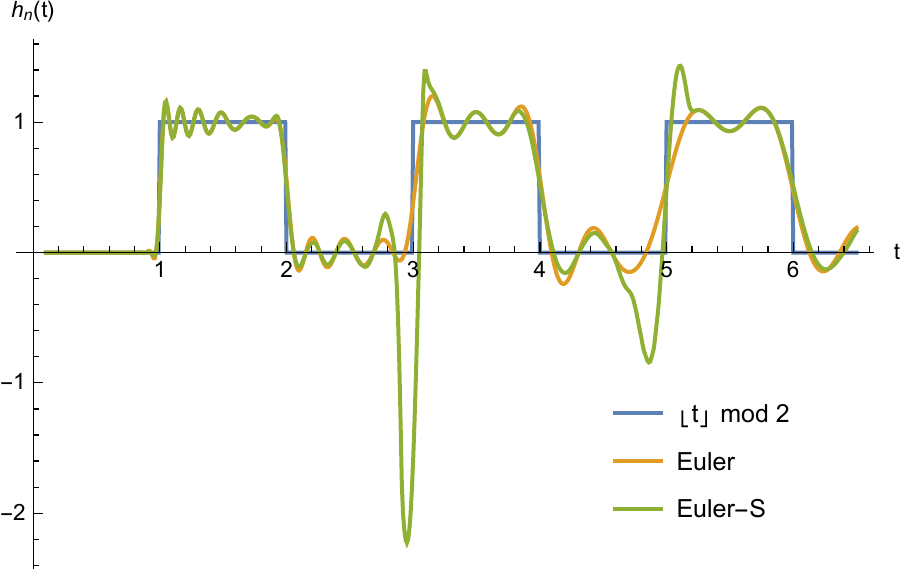}
\end{minipage}
\hfill
\caption{The Euler and the Euler-S NILT methods applied for $h(t)=\lfloor t \rfloor \!\!\mod 2$  with order $30$ }
\label{fig:eulers}
\end{figure*}

In the next section we are going to present several examples with decaying functions, where the Euler-S method provides accurate results, but we have to emphasize that the 
coincidence of the optimal $\hat{\theta}$ parameter for the CME and the Euler methods is not ensured in general. As an example Figure \ref{fig:eulers} demonstrate the risks of using the 
Euler-S method for ``regular'' functions for ``small'' $T$ (the associated $\hat{\theta}$ values are depicted in Figure \ref{fig:regfunc}). 

\section{Numerical analysis of the CME-S and the Euler-S methods}
\label{s:num2}

First, we check the behaviour of the CME-S and the Euler-S methods for the set of decaying functions in Table \ref{tab:test_functions}, because non-optimized NILT procedures easily fail in such cases. 

\setlength\extrarowheight{2pt}
\begin{table}[th]
	\centering
	\begin{tabular}{|c|c|c|}
		\hline
		$h(t)$ &  $h^*(s)$ & $a$\\
		\hline
		\hline
		$\exp(-t^2)$ & $\frac12 e^{(s/2)^2}\sqrt{\pi} ~\mathrm{Erfc}(s/2)$ & $-\infty$  \\
		\hline
		$\exp(-t)$ & $\frac{1}{1+s}$  & $-1$\\
		\hline
		$\exp(-\sqrt{t})$ & $\frac{1}{s}- \frac{s^{-3/2} e^{\frac{1}{4s}}}{2} \sqrt{\pi}
		~\mathrm{Erfc}\left(\frac{1}{2	\sqrt{s}}\right)$ & $0$ \\
		\hline
		$\frac{2}{(1+t)^3}$ & $1 - s - e^s s^2 ~\mathrm{Ei}(-s)$ & $0$ \\ 
		\hline
	\end{tabular}
	\caption{The set of decaying test functions, where  $\mathrm{Erfc}(z)=\frac{2}{\sqrt{\pi}} \int_{0}^{z} e^{-t^2} dt$ (error function) and  $\mathrm{Ei}(z)=\int_{-z}^{\infty} \frac{e^{-t}}{-t} dt$ (exponential integral function)}
	\label{tab:test_functions}
\end{table}

Table \ref{tab:shift_table} presents the results of the  CME, CME-S, Euler, and Euler-S procedures together with the theoretical value (``precise''), the number of NILT evaluations required for optimization of the shifting parameter (``iter.''), and the optimal value of the shifting parameter (``$\hat{\theta}$'') for order $30$. Based on Table \ref{tab:shift_table} we conclude that 
\begin{itemize}
\item non-optimized NILT easily fails to properly approximate the order of magnitude of functions decaying to zero, 
\item NILT with shifting provides a much better approximation for these decaying function, when the original method is inaccurate, 
\item the convex optimization procedure (with the stopping criteria $\theta_h-\theta_\ell<\epsilon=0.1$) terminates in $\approx20$ iterations, which means that the computational complexity of the shifting based NILT is $\approx20$ times higher than the one without shifting. 
\end{itemize}

\setlength\extrarowheight{1pt}
\begin{table*}[!th]
\centering
\begin{small}
\begin{tabular}{|@{\,}c@{\,}||@{\,}c@{\,}||@{\,}c@{\,}|@{\,}c@{\,}||@{\,}c@{\,}|@{\,}c@{\,}||@{\,}c@{\,}|@{\,}c@{\,}|}
\hline
order & precise & CME  & CME-S & Euler  & Euler-S & $\hat{\theta}$ & iter.   \\
\hline
\hline
\multicolumn{8}{|c|}{$h(t)=\exp(-t^2),T=5,a=-\infty$}\\
\hline
$30$ & $1.389E\!-\! 11$ & $8.739E\!-\! 6$ & $1.372E\!-\! 11$ & $-1.221E\!-\! 10$ & $1.389E\!-\! 11$  & $-49.94$ & $17$ \\
\hline
$60$ & $1.389E\!-\! 11$ & $1.356E\!-\! 6$ & $1.385E\!-\! 11$ &  $1.389E\!-\! 11$ & $1.389E\!-\! 11$ & $-49.96$ &  $17$  \\
\hline
\hline
\multicolumn{8}{|c|}{$h(t)=\exp(-t^2),T=10,a=-\infty$}\\
\hline
$30$ & $3.720E\!-\! 44$ & $5.515E\!-\! 6$ & $3.557E\!-\! 44$ & $3.889E\!-\! 9$ & $3.720E\!-\! 44$  & $-199.98$ & $20$ \\
\hline
$60$ & $3.720E\!-\! 44$ & $8.911E\!-\! 7$ & $3.681E\!-\! 44$ & $3.205E\!-\! 17$ & $3.720E\!-\! 44$  & $-199.95$ & $20$ \\
\hline
\hline
\multicolumn{8}{|c|}{$h(t)=\exp(-t),T=10,a=-1$}\\
\hline
$30$ & $4.540E\!-\! 5$ & $5.226E\!-\! 5$ & $4.540E\!-\! 5$ & $4.540E\!-\! 5$ & $4.540E\!-\! 5$ & $-10.01$ & $13$ \\
\hline
$60$ & $4.540E\!-\! 5$ & $4.654E\!-\! 5$ & $4.540E\!-\! 5$ & $4.540E\!-\! 5$ & $4.540E\!-\! 5$ & $-10.01$ & $14$ \\
\hline
\hline
\multicolumn{8}{|c|}{$h(t)=\exp(-t),T=50,a=-1$}\\
\hline
$30$ & $1.929E\!-\! 22$ & $2.111E\!-\! 6$ & $1.929E\!-\! 22$ & $-4.777E\!-\! 12$ & $1.929E\!-\! 22$ & $-49.99$ & $16$ \\
\hline
$60$ & $1.929E\!-\! 22$ & $3.273E\!-\! 7$ & $1.929E\!-\! 22$ & $-1.586E\!-\! 20$ & $1.929E\!-\! 22$ & $-49.99$ & $16$ \\
\hline
\hline
\multicolumn{8}{|c|}{$h(t)=\exp(-\sqrt{t}),T=100,a=0$}\\
\hline
$30$ & $4.540E\!-\! 5$ & $4.748E\!-\! 5$ & $4.544E\!-\! 5$ & $4.540E\!-\! 5$ & $4.540E\!-\! 5$ & $-5.35$ & $12$ \\
\hline
$60$ & $4.540E\!-\! 5$ & $4.573E\!-\! 5$ & $4.541E\!-\! 5$ & $4.540E\!-\! 5$ & $4.540E\!-\! 5$ & $-5.269$ & $12$ \\
\hline
\hline
\multicolumn{8}{|c|}{$h(t)=2/(1+t)^{3},T=100,a=0$}\\
\hline
$30$ & $1.941E\!-\! 6$ & $3.860E\!-\! 6$ & $1.954E\!-\! 6$ & $1.941E\!-\! 6$ & $1.941E\!-\! 6$ & $-5.72$ & $7$ \\
\hline
$60$ & $1.941E\!-\! 6$ & $2.214E\!-\! 6$ & $1.934E\!-\! 6$ & $1.941E\!-\! 6$ & $1.941E\!-\! 6$  & $-5.63$ &  $7$ \\
\hline
\end{tabular}
\end{small}
\caption{Properties of CME and Euler based NILT with and without shifting}
\label{tab:shift_table}
\end{table*} 

\begin{figure*}
\begin{minipage}{0.49\textwidth}
\centering\includegraphics[width=0.98\textwidth]{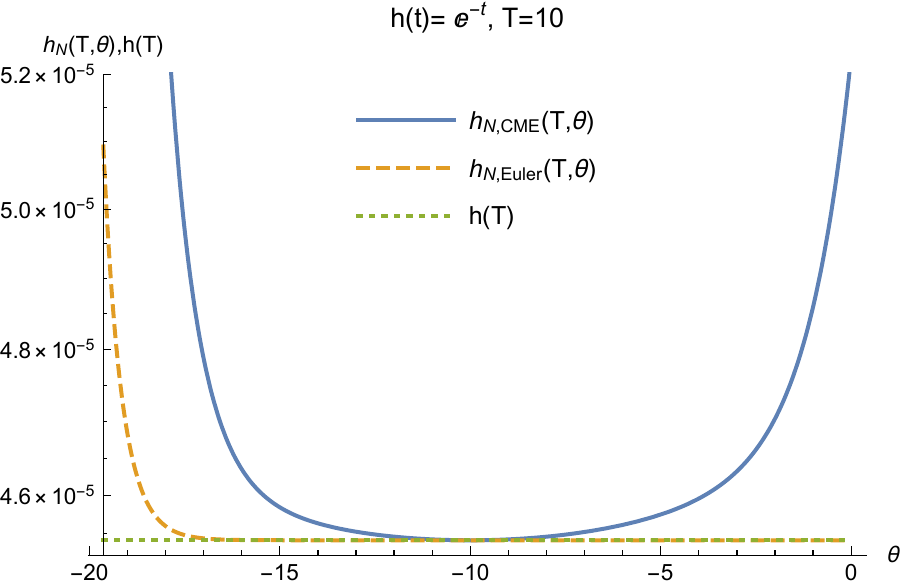}
\\ {\small a) $T=10$}
\end{minipage}
\hfill
\begin{minipage}{0.49\textwidth}
\centering\includegraphics[width=0.98\textwidth]{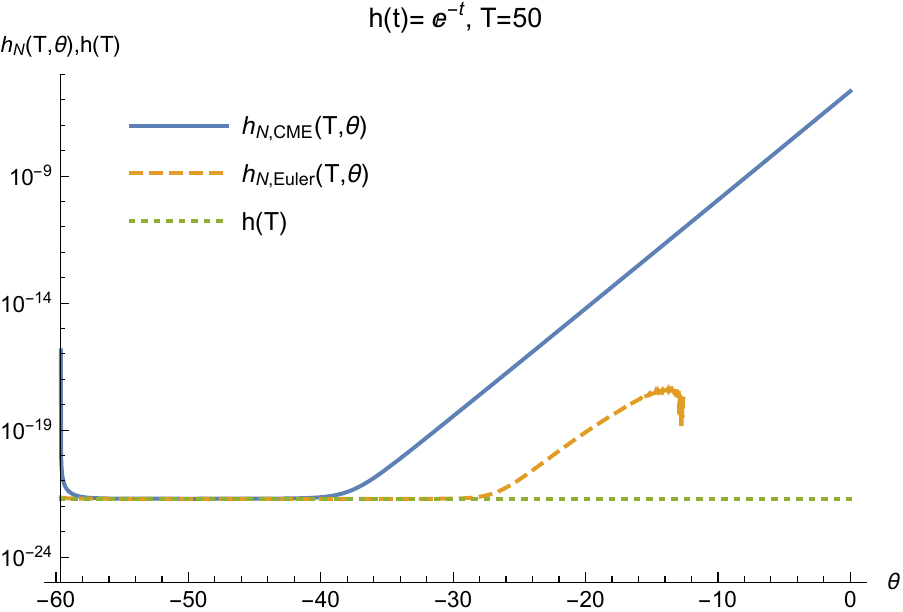}
\\ {\small b) $T=50$}
\end{minipage}
\caption{CME and Euler NILT with various shifting parameters when $h(t)=e^{-t}$, $T=10$  and $T=50$  with order $30$}
\label{fig:exptheta}
\end{figure*}

\begin{figure*}
\begin{minipage}{0.49\textwidth}
\centering\includegraphics[width=0.98\textwidth]{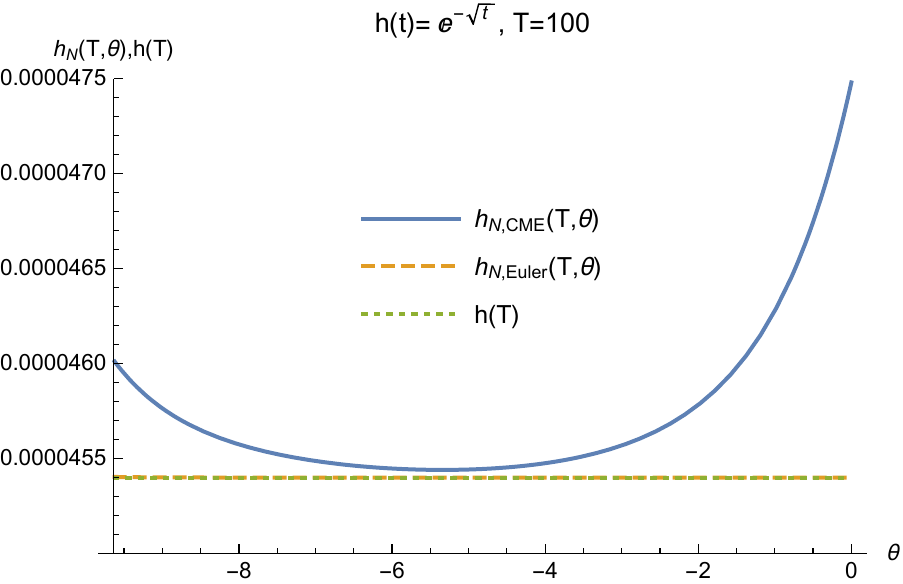}
\\ {\small a) $h(t)=e^{-\sqrt{t}}$}
\end{minipage}
\hfill
\begin{minipage}{0.49\textwidth}
\centering\includegraphics[width=0.98\textwidth]{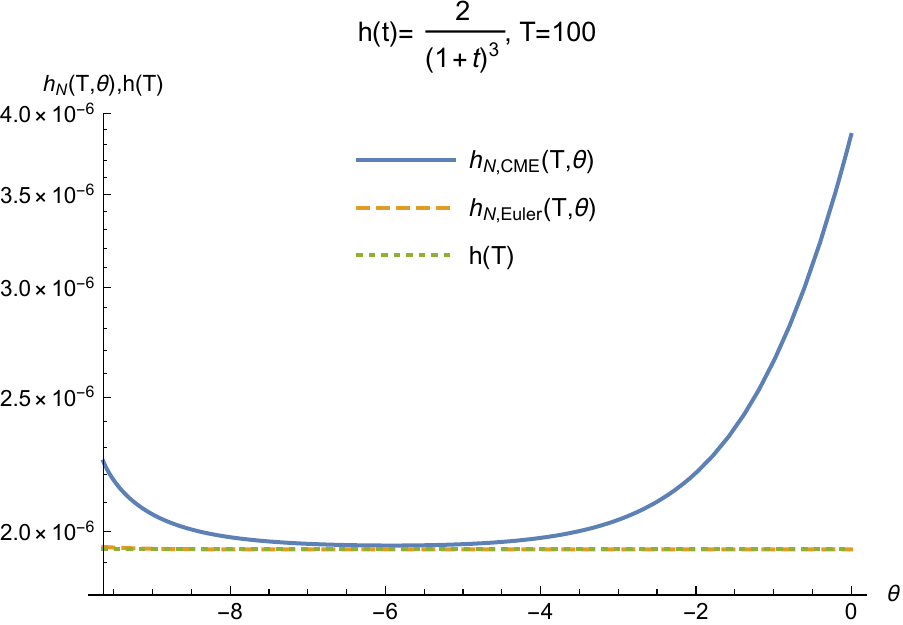}
\\ {\small b) $h(t)=\frac{2}{(1+t)^3}$}
\end{minipage}
\caption{CME and Euler NILT with various shifting parameters when  $T=100$,  $h(t)=e^{-\sqrt{t}}$ and $h(t)=\frac{2}{(1+t)^3}$ with order $30$}
\label{fig:gyokpartheta}
\end{figure*}

The effect of shifting for the other evaluated cases of Table \ref{tab:shift_table} is depicted in Figures  \ref{fig:exptheta} and \ref{fig:gyokpartheta}. 
In case of $h(t)=\exp(-t),T=10$, the accurate $\theta$ region of the Euler method is wide enough to contain $\theta=0$ and consequently the original Euler method provides accurate result. In case of $h(t)=\exp(-t),T=50$ the accurate $\theta$ region of the Euler method ends at $\theta=-30$ and the original Euler method provides a negative NILT result. 
At the optimized $\hat{\theta}$ value both methods are accurate (c.f. Figures  \ref{fig:exptheta}). 
In case of $h(t)=\exp(-\sqrt{t}),T=100$ and $h(t)=\frac{2}{(t+1)^3},T=100$, the accurate $\theta$ region of the Euler method is wide enough to contain $\theta=0$, i.e., the original Euler method is accurate enough, and the optimized CME method provides similarly accurate result (c.f. Figure \ref{fig:gyokpartheta}).

\section{Conclusion}
\label{s:concl}

In this paper we consider two efficient NILT methods of the AWF, the Euler and the CME method, and propose their enhancement with an optimized shifting parameter, which depends on the transform function and the time point of interest. The enhanced procedures are referred to as Euler-S ans CME-S. 

The paper presents many examples for the behaviour of these NILT methods, where the NILT results are compared with the (known) inverse Laplace values and provides intuitive explanations for their features.  
The goal of general purpose NILT is to provide trustable NILT approximate for any Laplace domain function without detailed knowledge on expected NILT result. For such cases we propose the use of the CME-S method which performs well in a wide range of the cases. 

If some background information is available about the expected behaviour of the inverse Laplace function, one can make a more accurate choice of the applied NILT method. 
In such cases our proposal is to use the Euler method for ``smooth'' functions with small $T$ (e.g., where $h(t)>10^{-10}$), use the Euler-S method for ``smooth'' functions with large with  $T$, and use the CME-S method otherwise. In any case handle the obtained result with special care if $a>-\infty$ and $\hat{\theta}=\theta_\ell$, because it might mean that the shifting parameter which balances the right and the left error is not feasible according to Assumption A3).

\bibliographystyle{ACM-Reference-Format}
\bibliography{references}

\end{document}